\documentclass{amsart}
\usepackage{amsmath}
\usepackage{amssymb}
\usepackage{amscd}
\usepackage{latexsym}
\usepackage{amsfonts}
 \newtheorem{thm}{Theorem}[section]
 \newtheorem{cor}[thm]{Corollary}
 \newtheorem{lem}[thm]{Lemma}
 \newtheorem{prop}[thm]{Proposition}
 \theoremstyle{definition}
 \newtheorem{defn}[thm]{Definition}
\theoremstyle{remark}
\newtheorem{rem}[thm]{Remark}
\newtheorem{ex}[thm]{Example}
\numberwithin{equation}{section}
\DeclareMathOperator{\ran}{ran}
\DeclareMathOperator{\spa}{span}
\begin{document}

%--------------------------------------------------------------------------
% editorial commands: to be inserted by the editorial office %
%\firstpage{1} %\issuenumber{4} %\Volumeandyear{41 (2001)}
%\Copyrightyear{2002} %\Signet %\commby{inhouse} %\submitted{March 14,
%2000} %\received{March 16, 2000} %\revised{June 1, 2000}
%\accepted{July 22, 2000}
%---------------------------------------------------------------------------
%Insert here the title, affiliations and abstract: %

\title{Rational
hyperholomorphic functions in \(\mathbb R^4\)}
\author[D. Alpay]{Daniel Alpay} \address{Department of Mathematics\\
Ben-Gurion University of the Negev\\ Beer-Sheva 84105\\ Israel}
\email{dany@math.bgu.ac.il} \thanks{The research of Daniel Alpay was
supported by the Israel Science Foundation (Grant no. 322/00).}
%-----Author 2

\author[M. Shapiro]{Michael Shapiro} \address{ Departamento de
Matem\'aticas\\ E.S.F.M. del I.P.N.\\ 07300 M\'exico, D.F.\\
M\'exico} \email{shapiro@esfm.ipn.mx} \thanks{The research of
Michael Shapiro was partially supported by CONACYT projects as
well as by Instituto Polit\'ecnico Nacional in the framework of
COFAA and CGPI programs}

\author[D. Volok]{Dan Volok} \address{Department of Mathematics\\
Ben-Gurion University of the Negev\\ Beer-Sheva 84105\\ Israel}
\email{volok@math.bgu.ac.il} 

\date{}

\begin{abstract} We introduce  the notion of rationality
 for hyperholomorphic functions (functions in the kernel of the
Cauchy-Fueter operator).  Following the case  of one complex
variable, we give three equivalent definitions: the first in terms
of Cauchy-Kovalevskaya quotients of polynomials, the second in
terms of realizations and the third in terms of backward-shift
invariance. Also introduced and studied are the counterparts of
the Arveson space and  Blaschke factors.
\end{abstract}
\maketitle
%\tableofcontents
%\setlength{\parskip}{1ex plus 0.5ex minus 0.2ex}
\section{Introduction} \label{intro} It is well known that
 functions holomorphic in a domain $\Omega\subset{\mathbb C}$ are exactly
the elements of the kernel of the Cauchy-Riemann differential operator
\begin{equation*} \overline{\partial}=\frac{\partial} {\partial
x}+i\frac{\partial}{\partial y} \end{equation*} restricted to
$\Omega$. A polynomial in $x$ and $y$ is holomorphic if,
 and only if, it is a polynomial in the complex variable $z=x+iy$, and
rational holomorphic functions are quotients of polynomials.\\

Holomorphic functions of one complex variable have a natural
generalization to the quaternionic setting when one replaces the
Cauchy-Riemann operator
 by the  Cauchy-Fueter operator \[D= \frac{\partial }{\partial x_0}+
{\mathbf e}_{\mathbf 1}\frac{\partial }{\partial x_1}+ {\mathbf
e}_{\mathbf 2}\frac{\partial }{\partial x_2}+ {\mathbf e}_{\mathbf
3}\frac{\partial }{\partial x_3}.  \] In this expression  the
$x_j$ are real variables and the ${\mathbf e}_{\mathbf j}$ are
imaginary units of the skew-field ${\mathbb H}$ of quaternions
(see Section \ref{federal triangle} below for more details).
Solutions of the equation $Df=0$ are called left-hyperholomorphic
functions (they are also called left-hyperanalytic, or
left-monogenic, or regular,  functions, see \cite{sprossig},
\cite{bds}, \cite{MR95m:30064}). Right-hyperholomorphic functions
are the solutions of the equation
\[fD= \frac{\partial f}{\partial x_0}+ \frac{\partial f}{\partial
x_1}{\mathbf e}_{\mathbf 1}+ \frac{\partial f}{\partial x_2}{\mathbf
e}_{\mathbf 2}+ \frac{\partial f}{\partial x_3}{\mathbf e}_{\mathbf
3}=0.  \]

When trying to generalize the notions of polynomial and rational
functions to the hyperholomorphic setting, one encounters several
obstructions. For instance, the quaternionic variable
\begin{equation*} x=x_0+x_1{\mathbf e}_{\mathbf 1}+x_2{\mathbf
e}_{\mathbf 2}+ x_3{\mathbf e}_{\mathbf 3} \end{equation*} is not
hyperholomorphic. Moreover, the point-wise product of two
hyperholomorphic functions is not hyperholomorphic in general and the
point-wise inverse of a non-vanishing hyperholomorphic function need not be
hyperholomorphic.\\

For the polynomials these difficulties were  overcome by Fueter, who
introduced  in \cite{fueter} the symmetrized multi-powers of the three
elementary functions \begin{equation*} \zeta_1(x)=x_1-{\bf
e_1}x_0,\quad\zeta_1(x)=x_2-{\bf e_2}x_0,\quad {\rm
and}\quad\zeta_3(x)=x_3-{\bf e_3}x_0.  \end{equation*} The polynomials
thus obtained are known today as the Fueter polynomials. They are
(both right and left) hyperholomorphic and appear in power series
expansions of hyperholomorphic functions.  In particular, a
hyperholomorphic polynomial is a linear combination of
the Fueter polynomials.\\

In this paper we introduce the notion of rational hyperholomorphic function.
 We obtain three equivalent characterizations: the
first one in terms of quotients and products of polynomials, the
second one in terms of realization and the last one in terms of
backward-shift-invariance. These various notions need to be
suitably defined in the hyperholomorphic setting. A key tool here
is the Cauchy-Kovalevskaya product of hyperholomorphic functions.
\\

We also
introduce a reproducing kernel Hilbert space of
left-hyperholomorphic functions which seems to be the
counterpart of the Arveson space of the ball --
 the reproducing kernel Hilbert space of functions
holomorphic in the open unit ball of ${\mathbb C}^N$  with the {\em rational}
reproducing
kernel $\frac{1}{1-\sum z_j\overline w_j}$ .
When $N=1$, this is just the Hardy space of the open unit disk.
 It was first
introduced by S. Drury in \cite{MR80c:47010} and
 proved in recent years to be a better extension of the
Hardy space than the classical Hardy space of the unit ball of
${\mathbb C}^N$, at least for problems in operator theory
 (see for instance
the papers \cite{agler-carthy}, \cite{adr1}, \cite{MR1758582},
 \cite{btv},  \cite{MR1903737}
for a sample of examples and applications).
In particular, it
 is
invariant under the operators $M_{z_j}$ of multiplication by the
variables $z_j,\,j=1,\ldots, N$, and it holds that
\begin{equation}
\label{2-novembre-2003}
I-\sum_1^N M_{z_j}M_{z_j}^*=C^*C
\end{equation}
where $C$ is the point evaluation at the origin.\\

To explain our approach let us consider briefly first the case of
 holomorphic  functions of one complex variable. Let $f$ and $g$ be two
functions holomorphic in a neighborhood of the origin, with the power
series expansions \begin{equation}\label{21-janvier-2003}
f(z)=\sum_{n=0}^\infty z^na_n\quad{\rm and}\quad
g(z)=\sum_{n=0}^\infty z^nb_n\end{equation} at the origin. Then the
point-wise product $(fg)(z)=f(z)g(z)$ has  at the origin the expansion
\begin{equation}\label{21-janvier-2003!}  (fg)(z)=\sum_{n=0}^\infty
z^nc_n, \end{equation} where
the sequence \(\{c_n\}\), given
by
\begin{equation}\label{21-janvier-2003!!}
c_n:=\sum_{m=0}^na_mb_{n-m},  \end{equation}
 is called the  convolution of the sequences \(\{a_n\}\) and \(\{b_n\}\).
It appears
that the substitute for pointwise product in the hyperholomorphic
setting (the Cauchy-Kovalevskaya product) is also a convolution.
\\

In the sixties of the previous century, the state space theory of
linear systems gave rise to a representation of a rational function called
{\sl realization} (see \cite{kfa}, \cite{bgk1}).
Still assuming analyticity in a neighborhood of the origin, this
representation is of the form \begin{equation}\label{yokbaba}
r(z)=D+zC(I-zA)^{-1}B \end{equation} where $A,B,C,D$ are matrices of
appropriate dimensions. It is particularly suitable for the study of
matrix-valued rational functions.\\

Realization theory has various extensions  in the setting of
several complex variables; see e.g. \cite{MR1855072},
\cite{MR2001e:93002}. One approach, related to  functions
holomorphic in the unit ball, exploits the so-called Gleason
problem (see \cite{akap1}, \cite{ad4}). A solution of the Gleason
problem, due to Leibenson (see \cite{MR87j:32012}, \cite[\S 15.8,
p.151]{MR54:11066}),
 was adapted to the setting of hyperholomorphic functions  in
\cite{as1} and \cite{as2}. It leads naturally to the analogues of
\eqref{21-janvier-2003} -- \eqref{21-janvier-2003!!}; these are the
expansions in terms of Fueter polynomials and the Cauchy-Kovalevskaya
product, mentioned above. Moreover, in this way we obtain the
analogue of the realization \eqref{yokbaba} and other equivalent
descriptions of the class of rational hyperholomorphic functions,
as well as the reproducing kernel of the counterpart of the  Arveson space
(quite different from the quaternionic Cauchy kernel). \\

This paper is organized as follows. In  Section \ref{federal
triangle} we review facts from the quaternionic analysis and
present some preliminary results, concerning backward-shift
operators in the hyperholomorphic setting.
 In Section \ref{smithsonian} we give three  definitions of a
rational function in the hyperholomorphic case and prove their
equivalence.
In  Section
\ref{hardy}
we define and study  the  counterparts
of the Arveson space of the unit ball
  and the Blaschke factors.\\

Some of the results presented here were announced in \cite{assv}.
In forthcoming papers we will consider the theory of linear
systems in the quaternionic case and Beurling-Lax-type theorems
for the Arveson space in the present setting.

\section{Quaternions and hyperholomorphic functions} \label{federal
triangle} \subsection{The skew-field of quaternions} \label{bassora}

In this section, we provide  some background on quaternionic analysis
needed in this paper. For more information, we refer the reader to
\cite{MR80g:30031} and to \cite{as3}.
The Hamilton skew-field of quaternions ${\mathbb H}$ is the real
four-dimensional linear space
 $\mathbb{R}^4$ equipped with the  product, defined
as follows.

  For the elements of the standard basis \({\bf e_0},{\bf e_1},{\bf
e_2},{\bf e_3}\)
 the rules of multiplication form the Cayley table:
\begin{equation}\label{cayley} \begin{array}{|c|c|c|c|c|} \hline
   &\bf e_0&\bf e_1&\bf e_2&\bf e_3\\ \hline
 \bf e_0&\bf e_0 &\bf e_1&\bf e_2&\bf e_3\\ \hline \bf e_1&\bf
e_1&-\bf e_0&\bf e_3&-\bf e_2\\ \hline \bf e_2&\bf e_2&-\bf
e_3&-\bf e_0&\bf e_1\\ \hline \bf e_3&\bf e_3&\bf e_2&-\bf
e_1&-\bf e_0 \\ \hline \end{array} \end{equation} Given two
elements \begin{align*} x=\sum^3_{i=0} x_i{\bf  e_i},&\  x_i\in
{\mathbb R},\\ y=\sum^3_{j=0} y_j{\bf  e_j},&\  y_j\in {\mathbb
R}, \end{align*}
 of ${\mathbb H},$ their  product is defined by
\[xy:=\sum^3_{i,j=0}x_iy_j{\bf e_ie_j},\] where \({\bf e_ie_j}\) are
calculated according to \eqref{cayley}.  Note that
 \(\bf e_0\) is the identity
element of  \(\mathbb{H}\) (for convenience, we identify it with the
real unit: \({\bf e_0}=1\)).

The quaternionic modulus \(|\cdot|\) coincides with the Euclidean norm
in $\mathbb{R}^4$:
\[|x|=\|x\|_{{\mathbb R}^4}=\sqrt{ \sum^3_{k=0} x_k^2},\]
and it holds that
\[|xy|=|x||y|\quad \forall x,y\in\mathbb H.\]

 The
conjugation in \(\mathbb{H}\) is defined by \[ \overline{x} =
 x_0 - \sum^3_{i=1} x_i {\bf  e_i}. \] It holds that
\[ \overline{x}x =
x\overline{x} =  |x|^2\] and hence
\[\forall x\in\mathbb H\setminus\{0\}\ :\
x^{-1}=\overline{x} |x|^{-2}.\]

\subsection{Hyperholomorphic functions and the Cauchy-Kovalevskaya
product} We have already mentioned in Section \ref{intro} that an
\(\mathbb{H}\)-valued function $f,$ \(\mathbb R\)-differentiable in an
 open connected set  $\Omega\subset\mathbb H,$ is  said to be
left-hyperholomorphic in \(\Omega\)  if it  satisfies in \(\Omega\)
the following differential equation: \begin{equation}\label{CFE}
\sum_{i=0}^3{\bf e_i}\frac{\partial f}{\partial x_i} = 0.
\end{equation} Analogously, an \(\mathbb{H}\)-valued function $f,$
\(\mathbb R\)-differentiable in an
 open connected set  $\Omega\subset\mathbb H,$ is  said to be
right-hyperholomorphic in \(\Omega\)  if it  satisfies in \(\Omega\)
the  differential equation \begin{equation}\label{CFEr}
\sum_{i=0}^3\frac{\partial f}{\partial x_i} {\bf e_i}= 0.
\end{equation}
The differential operator \[D=\sum_{i=0}^3{\bf e_i}\frac{\partial
}{\partial x_i}\] is called the Cauchy-Fueter operator. It
 satisfies the identity \[D\overline D=\overline
DD=\Delta_4,\] where
$$\overline D=\frac{\partial
}{\partial x_0}-\sum_{j=1}^3{\bf e_j}\frac{\partial }{\partial x_j}\quad{\rm
and}\quad
\Delta_4=\sum_{i=0}^3\frac{\partial^2}{\partial x_i^2}.$$
Thus hyperholomorphic functions are, in particular, harmonic.\\

In the sequel we shall restrict ourselves to
the case of left-hyperholomorphic functions. One can, of course,
obtain analogous results for
right-hyperholomorphic functions, as well.\\

   Let us denote the right-\(\mathbb H\)-module of functions,
left-hyperholomorphic in \(\Omega\), by \(\mathcal{O}_{\mathbb
H}(\Omega)\). Assume that \(\Omega\) is a ball, centered at the
origin. Then, as was proved in \cite{as2}, any element
\(f\in\mathcal{O}_{\mathbb H}(\Omega)\) can be written in the form
\begin{equation}\label{bientot-la-fin} f(x)=f(0)+
\sum_{n=1}^3\zeta_n(x)\mathcal R_nf(x), \end{equation} where
\begin{equation} \label{fueter-2003} \zeta_n(x):=x_n-x_0{\bf e_n}
\end{equation} are entire (both right and left) hyperholomorphic
functions, and the operators
\[\mathcal R_n:\mathcal O_{\mathbb H}(\Omega)\mapsto
\mathcal O_{\mathbb H}(\Omega)\]
 are defined by
 \begin{equation}
 \label{glesol}
 \mathcal R_nf(x)=\int_0^1\frac{\partial f}{\partial x_n}(tx)dt.
\end{equation}
(see \cite[p. 118]{rudin-ball}, \cite[\S 15.8 p.151]{MR54:11066}
and \cite{ad4} for these operators in the setting of the unit ball
of ${\mathbb C}^N$). Note that it follows from the
hyperholomorphic Cauchy integral formula that
\(\mathcal{O}_{\mathbb H}(\Omega) \subset C^{\infty}(\Omega)\),
hence the operators \(\mathcal R_n\) commute:
\[\mathcal R_m\mathcal R_n f(x)=\int_0^1\int_0^1
\frac{\partial^2 f}{\partial x_n\partial x_m}(utx)tdtdu
=\mathcal R_n\mathcal R_m f(x),\]
and that
\[ \mathcal R_nf(0)=\frac{\partial f}{\partial x_n}(0).\]
Hence,
applying the formula \eqref{bientot-la-fin} for \(\mathcal R_nf\),
we get
\[f(x)=f(0)+
\sum_{n=1}^3\zeta_n(x)\frac{\partial f}{\partial x_n}(0)+
\sum_{0\leq n\leq m\leq 3}(\zeta_n(x)\zeta_m(x)+\zeta_m(x)\zeta_n(x))
\mathcal R_m\mathcal R_n f(x).\]
Iterating this process,  one obtains an
expansion of \(f\) in terms of  symmetrized products of $\zeta_n,$
analogous to the classical Taylor power series expansion.
\\

To be more precise, let us introduce the multi-index notation
we shall use throughout this paper.
The symmetrized
 product  of  $a_1, \ldots, a_n\in\mathbb{H}$
is defined by
\[ a_1\times a_2\times\cdots\times a_n = \frac{1}{n!}\sum_{\sigma\in
S_n} a_{\sigma(1)}a_{\sigma(2)}\cdots a_{\sigma(n)}, \] where
$S_n$ is the set of all permutations of the set $\{1,...,n\}$.
Furthermore, for \(\nu,\mu\in\mathbb{Z}_+^3\) we use the usual
notation
\[|\nu|=\nu_1+\nu_2+\nu_3,\quad \nu! = \nu_1!\nu_2!\nu_3!, \quad
\nu\geq\mu\text{ if }\nu_j\geq\mu_j\ \forall j,\]
\[{\partial}^{\nu} =\frac{{\partial}^{|\nu|} } {\partial x^{\nu_1}_1
\partial x^{\nu_2}_2 \partial x^{\nu_3}_3},\]
 \[
 e_1=\begin{pmatrix}1&0&0\end{pmatrix},\quad
 e_2=\begin{pmatrix}0&1&0\end{pmatrix},\quad
 e_3=\begin{pmatrix}0&0&1\end{pmatrix}.
 \]

Using the above notation, we can formally write
\begin{equation}\label{uuab} f(x) = \sum_{n=
0}^\infty\sum_{|\nu|=n}\zeta^\nu(x)  f_\nu, \end{equation} where
\begin{align}\label{Fupol} \zeta^{\nu}(x)&:=
\zeta_1(x)^{\times\nu_1}\times
 \zeta_2(x)^{\times\nu_2}\times\zeta_3(x)^{\times\nu_3},\\ \label{Taco}
f_\nu &:= \frac{1}{\nu!}({\partial}^{\nu} f)(0).  \end{align} The
polynomials \(\zeta^\nu\), defined by \eqref{Fupol}, are called
the Fueter polynomials. It can be proved that the Fueter
polynomials are entire (both left and right) hyperholomorphic and
that the series \eqref{uuab} is normally convergent. Thus one can
characterize the right-\(\mathbb{H}\)-module
\(\mathcal{O}_{\mathbb{H}}\) of functions, left-hyperholomorphic
in a neighborhood of the origin, as follows (see \cite{bds}):

\begin{thm}\label{Taylor} An \(\mathbb{H}\)-valued function $f,$
defined in a neighborhood of the origin, belongs to the space
\(\mathcal{O}_{\mathbb{H}}\) if, and only if, it can be represented in
the form \eqref{uuab}, where \begin{equation}\label{ROC}
\rho(f)=\limsup_{n\to\infty}\left(\sum_{|\nu|=n}|f_\nu|\right)
^{\frac{1}{n}}<\infty. \end{equation}
In this case the series \eqref{uuab} converges uniformly
on compact subsets of the ball
\[\{x\in\mathbb{H}:|x|\cdot\rho(f)<1\}.\]
\end{thm} \begin{cor} An
\(\mathbb H\)-valued polynomial \(p\) of real variables
\(x_0,x_1,x_2,x_3\) is left-hyperholomorphic if, and only if, it is a
finite linear combination of Fueter polynomials:
\[p(x)=\sum_{n=0}^m\sum_{|\nu|=n}\zeta^{\nu}p_\nu,\quad p_\nu\in\mathbb H.\] \end{cor}

 \begin{rem}
 \label{analog}
In view of Theorem \ref{Taylor}, in the quaternionic analysis the
elementary functions \(\zeta_n\) play role, similar in a sense to
that of \(z_n\) in several complex variables. Thus \(\zeta_n\) are
sometimes called the hyperholomorphic variables. The term "total
variables" is used also referring to the fact that both
\(\zeta_n\) and all its powers are hyperholomorphic, see
\cite{bds}, \cite{sprossig}. We note, however, that \(\zeta_n\)
are neither independent, nor \(\mathbb H\)-linear. Moreover, the
choice of {\em left}-hyperholomorphic variables is not unique:
e.g., \(\zeta_n(\mathbf{e_1}x)\) are also suitable for this role,
but are not right-hyperholomorphic.
\end{rem}

It is  useful to  calculate the expressions for the operators
\(\mathcal R_n\), defined by \eqref{glesol},
  in terms of expansions \eqref{uuab}.

\begin{lem}\label{bsoplem} Let \(f\in\mathcal O_{\mathbb H}(\Omega)\)
be given by \eqref{uuab}. Then \begin{equation} \label{bsop} \mathcal
R_n f(x)=\sum_{\nu\geq e_n}\frac{\nu_n}{|\nu|}\zeta^{\nu-e_n}(x)
f_\nu.  \end{equation} \end{lem}

\begin{proof} Without loss of generality, \(f\) is a Fueter
polynomial. But \[\frac{\partial \zeta^\nu}{\partial
x_n}(x)=\nu_n\zeta^{\nu-e_n}(x),\] hence \[ \mathcal R_n
\zeta^\nu(x)=\int_0^1\frac{\partial \zeta^\nu}{\partial x_n}(tx)dt
=\int_0^1 \nu_n
t^{|\nu|-1}\zeta^{\nu-e_n}(x)dt=\frac{\nu_n}{|\nu|}\zeta^{\nu-e_n}(x).
\] \end{proof}

In view of Lemma \ref{bsoplem}, we propose the following \begin{defn}
\label{bsdef}
 The operators \(\mathcal R_n:\mathcal O_{\mathbb H}\mapsto\mathcal
O_{\mathbb H}\),
 defined by \eqref{bsop}, are called the backward-shift operators.
\end{defn}

 Following the analogy with the complex case, we would like to impose
on \(\mathcal{O}_{\mathbb{H}}\) the structure of a ring. However, the
point-wise product is not suitable here. For instance, the function
\(\zeta_1\zeta_2\) is not hyperholomorphic. Instead, one can use  (see
\cite[Section 14]{bds} and compare with \eqref{21-janvier-2003} --
\eqref{21-janvier-2003!!} in Section \ref{intro}) the following
\begin{defn} The Cauchy-Kovalevskaya product (below: C-K-product)
\(f\odot g\) of the functions \[f=\sum\zeta^\nu f_\nu,\
g=\sum\zeta^\nu g_\nu, \] left hyperholomorphic in a neighborhood
of the origin, is defined by
\begin{equation}\label{abc} f\ {\odot}\ g =
\sum_{n=0}^\infty\sum_{|\eta|=n}\zeta^\eta \sum_{0\leq \nu \leq
\eta} f_\nu g_{\eta-\nu}.  \end{equation} \end{defn}

\begin{rem}
In certain special cases the C-K-product coincides with the
point-wise one. For instance, if \(g(x)\equiv\text{const}\) then
\(f\odot g=fg\), but not necessarily \(g\odot f=gf\)! Another
special case is discussed in Section \ref{generali}.
\end{rem}

\begin{prop}\label{RON} The space \(\mathcal{O}_{\mathbb{H}}\),
equipped with the  C-K-product, is a ring. Moreover,
\[\rho(f\odot g)\leq\max\{\rho(f),\rho(g)\}.\] \end{prop} \begin{proof}
 Without loss of generality, we take \(\rho(f)=\rho(g)=\rho.\) Then \(
\forall\epsilon>0, \exists C(\epsilon)>0, \forall k: \)
\begin{align*}
 \sum_{|\nu|=k}|a_\nu|&\leq C(\epsilon)(\rho+\epsilon)^{k}, \\
\sum_{|\nu|=k}|b_\nu|&\leq  C(\epsilon)(\rho+\epsilon)^{k}.
\end{align*} Hence
\begin{multline*}\sum_{|\eta|=n}
|\sum_{0\leq \nu \leq \eta} a_\nu b_{\eta-\nu}|\leq
\sum_{|\eta|=n}\sum_{0\leq \nu \leq \eta} |a_{\nu}||b_{\eta-\nu}|
\leq\sum_{k=0}^n\sum_{|\nu|=k}|a_\nu|\sum_{|\mu|=n-k}|b_{\mu}|\\
\leq(n+1)C(\epsilon)^2(\rho+\epsilon)^n \end{multline*} and so
$$\rho(f\odot g)\leq\rho.$$ \end{proof} The C-K-product can be
generalized to spaces of matrix-valued left-hyperholo\-morphic
functions in the usual way: for \[F=(f_{\alpha,\beta})\in
\mathcal{O}^{m\times n}_{\mathbb{H}},\ G=(g_{\beta,\gamma})\in
\mathcal{O}^{n\times p}_{\mathbb{H}}\] we define \[F\odot
G:=\left(\sum_{\beta}f_{\alpha,\beta} \odot
g_{\beta,\gamma}\right)_{\alpha,\gamma}.\] The question arises,
when an element \(F\in\mathcal{O}^{n\times n}_{\mathbb{H}}\) is
C-K-invertible. In view of \eqref{abc},
 a necessary condition is that
the value \(F(0)\) must be invertible in \(\mathbb H^{n\times
n}\). This turns out to be also sufficient: \begin{prop}
\label{damascus1} Let \(F\in\mathcal{O}^{n\times n}_\mathbb{H}\).
If \(F(0)\) is invertible in \(\mathbb H\) then \(F\) is
C-K-invertible in \(\mathcal{O}^{n\times n}_\mathbb{H}\) and its
C-K-inverse $F^{-\odot}\in\mathcal{O}^{n\times n}_\mathbb{H}$  is
given by the series
\begin{equation}
F^{-\odot}=(F(0))^{-1}\odot(I_n-G)^{-\odot}=(F(0))^{-1}\odot
\sum_{k=0}^{\infty}G^{\odot k}, \end{equation} where
\[G=I_n-F(F(0))^{-1}.\] \end{prop}
\begin{proof}
 It suffices to show the normal convergence of the series
\begin{equation}\label{i-g}
(I_n-G)^{-\odot}=\sum_{k=0}^{\infty}G^{\odot k} \end{equation} in
a neighborhood of the origin for arbitrary
\(G\in\mathcal{O}^{n\times n}_\mathbb{H}\),\ satisfying
\(G(0)=0.\) According to Theorem \ref{Taylor}, \[
G=\sum_{p=1}^\infty\sum_{|\nu|=p}\zeta^\nu A_\nu,\] and there
exists \(A\in\mathbb{R}^+\) such that
\[\forall p>0:\quad\sum_{|\nu|=p}\|A_\nu\|\leq
A^p,\] where \(\|\cdot\|\)  denotes the operator norm.
Then \[ G^{\odot
k}=\sum_{p=k}^{\infty}\sum_{|\nu|=p}\zeta^{\nu}
\sum_{\substack{\mu_1+\ldots +\mu_k=\nu \\ \mu_1,\ldots,\mu_k\not=0}}
A_{\mu_1}\cdots A_{\mu_k}.  \] But for \(p\geq k\geq 1\)
\begin{eqnarray*} \sum_{|\nu|=p}\| \sum_{\substack{\mu_1+\ldots
+\mu_k=\nu \\ \mu_1,\ldots,\mu_k\not=0}} A_{\mu_1}\cdots
A_{\mu_k}\|&\leq&\sum_{\substack{p_1+\ldots+p_k=p\\
p_1,\ldots,p_k>0}}\ \sum_{|\mu_1|=p_1}\|A_{\mu_1}\|\cdots
\sum_{|\mu_k|=p_k}\|A_{\mu_k}\|\\ &\leq&\binom{p-1}{k-1}A^p<(2A)^p.
\end{eqnarray*} Therefore, \[|x|<\dfrac{1}{4A}\implies\|G^{\odot
k}(x)\|\leq\dfrac{1}{2^{k-1}},\] and the  normal convergence of the
series \eqref{i-g} in the ball \(\{x:|x|<1/4A\}\) follows.  \end{proof}
\subsection{The Gleason problem in the hyperholomorphic case}
\label{mossul}
In view of Remark \ref{analog},  the formula \eqref{bientot-la-fin}
may be considered as a solution for a Gleason problem with respect to
the hyperholomorphic variables \(\zeta_n\) (see \cite{as1}, \cite{as2}
for details and references).  However, there is a disadvantage in that
the point-wise product  appears. In particular,  the individual terms
\(\zeta_n(x)\mathcal R_nf(x)\) in the sum \eqref{bientot-la-fin} need
not be left-hyperholomorphic, in general. The goal of the present
section is to consider the Gleason problem with the point-wise product
being replaced by the C-K  product.

\begin{defn} \label{defgle}
 Let \(f\in\mathcal O_{\mathbb H}\).  The Gleason problem for \(f\) is
to find
 a triple of functions \(g_1,g_2,g_3\in\mathcal O_{\mathbb H}\), such
that \[f-f(0)=\sum_{n=1}^3\zeta_n\odot g_n.\] \end{defn}

It turns out that the backward-shift operators \(\mathcal R_n\)
provide a solution for this new Gleason problem, as well.

\begin{thm} \label{To kill a mockingbird} Let
$f\in\mathcal{O}_{\mathbb{H}}.$ Then it holds that \begin{equation}
\label{to kill a mockingbird} f-f(0)=\sum_{n=1}^3\zeta_n\odot \mathcal
R_nf.  \end{equation} \end{thm} \begin{proof} According to
\eqref{bsop}, we have \[\sum_{n=1}^3\zeta_n\odot \mathcal
R_nf=\sum_{n=1}^3\sum_{\nu\geq e_n}\frac{\nu_n}{|\nu|}\zeta^{\nu}
f_\nu=\sum_{|\nu|>0}\zeta^{\nu} f_\nu=f-f(0).\] \end{proof}

In general, the solution for the Gleason problem, provided  by the
backward-shift operators, is not the only possible one. To
illustrate this observation, let us consider the subspaces of
\(\mathcal{O}_{\mathbb{H}}^m\), in which the problem is solvable.

\begin{defn} \label{by Harper Lee}
 A subspace \(\mathcal{W}\) of \(\mathcal{O}^{m}_{\mathbb{H}}\) is said
to  be resolvent-invariant if \[\forall f\in\mathcal{W}\ \exists
g_1,g_2,g_3\in\mathcal{W}\ :\ f-f(0)=\sum_{n=1}^3\zeta_n\odot
g_n.\]
 If, moreover,  the space \(\mathcal{W}\)  is $\mathcal
R_n$-invariant for $n=1,2,3$, it is said to be
backward-shift-invariant.  \end{defn}

\begin{thm} \label{findiminv} A finite-dimensional  subspace
\(\mathcal{W}\) of \(\mathcal{O}^{m}_{\mathbb{H}}\) is
resolvent-invariant (respectively, backward-shift-invariant) if,
and only if, it is spanned by the columns of a matrix-valued
function of the form \begin{equation} \label{matker}
W=C\odot(I-\sum_{n=1}^3\zeta_nA_n)^{-\odot}, \end{equation}
 where \(C\) and \(A_n\) are constant matrices with entries in
${\mathbb H}$ (respectively, \(A_n\) commute).  \end{thm} In the
proof of Theorem \ref{findiminv} we shall use the following
\begin{lem}\label{kernel} Let $A_1,A_2$ and $A_3$ be in ${\mathbb
H}^{\ell\times \ell}$. Then in a neighborhood of the origin it
holds that \begin{equation}
(I_\ell-\zeta_1A_1-\zeta_2A_2-\zeta_3A_3)^{-\odot}=\sum_{\nu\in{\mathbb
Z}^3_+} \zeta^\nu A^\nu\frac{|\nu|!}{\nu!},\end{equation} where
\begin{equation}\label{denota}
A^\nu=A_1^{\times \nu_1}\times A_2^{\times \nu_2}\times
A_3^{\times \nu_3}.
\end{equation}
\end{lem} \begin{proof}
 We have
\[(I_\ell-\zeta_1A_1-\zeta_2A_2-\zeta_3A_3)^{-\odot}=\sum_{k=0}^\infty
\left(\zeta_1A_1+\zeta_2A_2+\zeta_3A_3\right)^{\odot k}.\] Let us
prove by induction on \(k\) that \begin{equation} \label{newton}
\left(\zeta_1A_1+\zeta_2A_2+\zeta_3A_3\right)^{\odot k}=
\sum_{|\nu|=k} \zeta^\nu\left(A_1^{\times \nu_1}\times A_2^{\times
\nu_2}\times A_3^{\times
\nu_3}\right)\frac{|\nu|!}{\nu!}.\end{equation} Indeed,
\eqref{newton} obviously holds for \(k=0,\) and if it holds for
some \(k\) then we have \begin{multline*}
\left(\zeta_1A_1+\zeta_2A_2+\zeta_3A_3\right)^{\odot(k+1)}
=\left(\zeta_1A_1+\zeta_2A_2+\zeta_3A_3\right)\odot\sum_{|\nu|=k}
\zeta^\nu A^\nu\frac{|\nu|!}{\nu!}\\
=\sum_{|\nu|=k+1}\zeta^\nu \frac{(|\nu|-1)!}{\nu!}\left(\nu_1A_1
A^{\nu-e_1} +\nu_2A_2A^{\nu-e_2}+\nu_3A_3A^{\nu-e_3}\right)\\
=\sum_{|\nu|=k+1}\zeta^\nu A^\nu\frac{|\nu|!}{\nu!}.
\end{multline*} \end{proof}

 \begin{proof}[Proof of Theorem
\ref{findiminv}] Let ${\mathcal W}$ be a resolvent-invariant
finite-dimensional
 subspace of \(\mathcal{O}^{m}_{\mathbb{H}}\) and let $W$ be a
matrix-valued hyperholomorphic function whose columns form a basis
of ${\mathcal W}$. Then there exist constant matrices
\(A_n\in{\mathbb H}^{\ell\times \ell}\) (with \(\ell=\dim{\mathcal
W}\)) and \(C=W(0)\), such that \[W=C+\sum_{n=1}^3\zeta_n\odot
WA_n=C+W\odot\sum_{n=1}^3\zeta_nA_n,\] hence \(W\)  is of the form
\eqref{matker}. If ${\mathcal W}$ is, moreover,
backward-shift-invariant then \(A_n\) can be chosen such that
$\mathcal R_nW=WA_n.$ Then $A_n$ commute since the $R_n$ do.

Conversely, let ${\mathcal W}$ be the span of the columns of a
matrix-valued function of the form \eqref{matker}. Then
\begin{multline*} W-W(0)=C\odot(I-\sum_{n=1}^3\zeta_nA_n)^{-\odot}-C\\
=W
\odot\left(I-(I-\sum_{n=1}^3\zeta_nA_n)\right)=\sum_{n=1}^3\zeta_n\odot
WA_n,
\end{multline*}
 and hence ${\mathcal W}$ is resolvent-invariant. If, moreover,
the matrices \(A_n\) commute then, according to Lemma \ref{kernel},
 \[W=\sum_{\nu\in\mathbb Z_+^3} \zeta^\nu CA_1^{\nu_1}A_2^{\nu_2}A_3^{
\nu_3}\frac{|\nu|!}{\nu!},\] hence $\mathcal R_nW=WA_n.$ This
completes the proof.  \end{proof}

\section{Rational hyperholomorphic functions} \label{smithsonian}
\subsection{Definitions} \label{kurdistan} In this section we give
three definitions of a rational function, left-hyperholomorphic in a
neighborhood of the origin. We prove that they are equivalent in
Section \ref{jean-yanne1}.\\

The first definition parallels the classical definition in terms
of quotients of polynomials in the complex case. Here polynomials
are replaced by the Fueter polynomials, point-wise multiplication
is replaced by the C-K-product, and inverses are replaced by the
C-K-inverses.

\begin{defn}\label{th1111} An \(\mathbb H^{m\times n}\)-valued
function \(R\), left-hyperholomorphic in a neighborhood of the
origin, is said to be rational if all its entries belong to the
minimal subring \(\mathcal Q_{\mathbb{H}}\) of
\(\mathcal{O}_{\mathbb{H}}\), which contains  hyperholomorphic
polynomials and is  closed under C-K-inversion: \[ r\in\mathcal
Q_{\mathbb{H}},r(0)\not=0 \implies \exists r^{-\odot}\in\mathcal
Q_{\mathbb{H}}.\] \end{defn}

\begin{ex}\label{Feb17!} Let $j=1,2,3$. The functions $\zeta_j$,
$\zeta^2_j$ and more generally all the Fueter polynomials are
rational.  \end{ex}

\begin{ex} The function $$ \left(\left((1-\zeta_1{\bf e}_{\mathbf
1})^{-\odot}+2\right)^{-\odot}+
\zeta_1\odot\zeta_2^{{\odot}3}\right)^{-\odot}+\zeta_3^{\odot
5}{\mathbf e}_{\mathbf 3}$$ is rational.\end{ex}

The next example will play an important role in the sequel.

\begin{ex} Let $a\in{\mathbb H}$. The function \begin{equation}
x\mapsto(1-\zeta_1\overline{\zeta_1(a)}-\zeta_2\overline{\zeta_2(a)}
-\zeta_3\overline{\zeta_3(a)})^{-\odot} \end{equation} is rational.
\end{ex}

The second definition parallels the realization \eqref{yokbaba} (see
Section \ref{intro}) in the complex case.  \begin{defn}\label{mum} An
\(\mathbb H^{m\times n}\)-valued function \(R\), left-hyperholomorphic
in a neighborhood of the origin, is said to be rational  if it can be
represented in the form \begin{equation}\label{voltaire}
R=D+C\odot(I-\zeta_1A_1-\zeta_2A_2-\zeta_3A_3)^{-\odot}\odot
(\zeta_1B_1+\zeta_2B_2+\zeta_3B_3), \end{equation} where $A_i,B_i$,
$C$ and $D$ are constant matrices with entries in ${\mathbb H}$ and of
appropriate dimensions.  \end{defn}

For brevity, from now on we shall use the  notation
 \begin{equation} \zeta_{(\ell)} :=\begin{pmatrix}\zeta_1 I_\ell&\zeta_2
I_\ell&\zeta_3 I_\ell \end{pmatrix}\in\mathcal O_{\mathbb
H}^{\ell\times 3\ell}.
\end{equation} The dimension $\ell$ will usually be understood
from the context and omitted. Then \eqref{voltaire} can be
rewritten as:
\begin{equation}\label{voltaire1} R=D+C \odot (I-\zeta A)^{-\odot} \odot
\zeta B \end{equation} where \begin{equation}
A=\begin{pmatrix}A_1\\A_2\\A_3\end{pmatrix}{\rm and}\quad
B=\begin{pmatrix}B_1\\B_2\\B_3\end{pmatrix}.  \end{equation}

 The
third definition is in terms of the resolvent-invariance.

\begin{defn}\label{defrat3} An \(\mathbb H^{m\times n}\)-valued
function \(R\), left-hyperholomorphic in a neighborhood of the origin,
is said to be rational  if there is a finite-dimensional
resolvent-invariant space \(\mathcal{W}\subset\mathcal{O}^{m}_{\mathbb
H}\),
 such that for every $v\in\mathbb{H}^n$ the Gleason problem for \(Rv\)
is solvable in $\mathcal{W}.$ \end{defn} The main result of the paper,
presented in Section \ref{jean-yanne1}, is that all three definitions
are equivalent. In view of Proposition \ref{kernel}, they are also
equivalent to the following

\begin{defn} \label{def4-prime} An \(\mathbb H^{m\times n}\)-valued
function \(R\), left-hyperholomorphic in a neighborhood of the origin,
is said to be rational
 if it can be represented as
\[R=\sum_{n=0}^\infty\sum_{|\nu|=n}\zeta^\nu R_\nu,\] where for
\(|\nu|\geq 1\) \[
R_\nu=\dfrac{(|\nu|-1)!}{\nu!}C\begin{pmatrix}\nu_1A^{\nu-e_1} &
\nu_2 A^{\nu-e_2} & \nu _3 A^{\nu-e_3}\end{pmatrix}B \]
with
\(A,B,C\) being constant matrices of appropriate dimensions.
\end{defn}

\subsection{Preparatory lemmas} \label{tchirchik}
The proof of the equivalence of Definitions \ref{th1111} -- \ref{defrat3}
is based on several technical lemmas.

\begin{lem}\label{nato} Let \(R\in\mathcal{O}^{n\times
n}_{\mathbb{H}}\) admit the representation \eqref{voltaire1},
where \(D\in\mathbb{H}^{n\times n}\) is invertible. Then \(R\) is
C-K-invertible and its C-K-inverse \(R^{-\odot}\) admits the
representation
  \begin{equation} R^{-\odot}=D^{-1}-D^{-1}C \odot
(I-Z\tilde A)^{-\odot} \odot ZBD^{-1}, \end{equation} where
$\tilde A=A-BD^{-1}C$.  \end{lem}
\begin{proof}
 We have: \begin{multline*}  (D+C\odot(I-\zeta A)^{-\odot}\odot
\zeta B) \odot (D^{-1} - D^{-1}C \odot (I-\zeta\tilde A)^{-\odot}
\odot \zeta BD^{-1}) \\ = I- C \odot(I-\zeta \tilde
A)^{-\odot}\odot \zeta BD^{-1} + C \odot (I-\zeta A)^{-\odot}
\odot \zeta BD^{-1} -\\ - C \odot(I-\zeta A)^{-\odot} \odot \zeta
BD^{-1}C \odot (I-\zeta \tilde A)^{-\odot}
 \odot \zeta BD^{-1} \\ = I - C \odot \big\{(I-\zeta \tilde A)^{-\odot} -
(I-\zeta A)^{-\odot} +\\ +(I-\zeta A)^{-\odot} \odot \zeta
BD^{-1}C \odot (I-\zeta \tilde A)^{-\odot}\big\} \odot \zeta
BD^{-1}.
\end{multline*} But
\[\zeta BD^{-1}C=\zeta (A-\tilde A)=(I-\zeta \tilde A)-(I-\zeta A),\] hence the expression
in the curly brackets is equal to \(0.\) \end{proof}

 \begin{lem}\label{natii} There exists a unitary
matrix $U\in{\mathbb H}^{3(\ell+ m)\times 3(\ell+ m)}$ such that
\begin{equation} {\rm diag}~(\zeta_{(\ell)},\zeta_{(m)})= \zeta_{(\ell+m)}U.
\end{equation} \end{lem}

\begin{proof} It  suffices to take $$
U=\begin{pmatrix} I_\ell&0&0&0&0&0\\ 0&0&0&I_m&0&0\\
0&I_\ell&0&0&0&0\\ 0&0&0&0&I_m&0\\ 0&0&I_\ell&0&0&0\\
0&0&0&0&0&I_m\end{pmatrix}. $$
\end{proof}
%\mbox{}\qed\mbox{}\\

\begin{lem}\label{nati} Let $R_i\in\mathcal{O}^{m_i\times
n_i}_{\mathbb{H}} $, $ i=1,2$, admit the representations
$$R_i(x)=D^{(i)}+C^{(i)} \odot (I-\zeta A^{(i)})^{-\odot} \odot
\zeta B^{(i)}.$$

If \(n_2=m_1\) then
\(R_1\odot R_2\) admits the representation
\begin{multline*}
 R_1 \odot R_2=D^{(1)}D^{(2)}+\\+\begin{pmatrix}
C^{(1)}&D^{(1)}C^{(2)}\end{pmatrix}\odot \left(I-\zeta
U\begin{pmatrix}A^{(1)}&B^{(1)}C^{(2)}\\0&A^{(2)}\end{pmatrix}
\right)^{-\odot} \odot \zeta
U\begin{pmatrix}B^{(1)}D^{(2)}\\B^{(2)}\end{pmatrix}.
\end{multline*}

If \(m_1=m_2\), \(n_1=n_2\) then \(R_1+ R_2\) admits the representation
\begin{multline*}
R_1+R_2=D^{(1)}+D^{(2)}+\\+\begin{pmatrix}
C^{(1)}&C^{(2)}\end{pmatrix}\odot \left(I-\zeta
U\begin{pmatrix}A^{(1)}&0\\0&A^{(2)}\end{pmatrix}
\right)^{-\odot}\odot \zeta
U\begin{pmatrix}B^{(1)}\\B^{(2)}\end{pmatrix}.
\end{multline*}
In both formulas $U$ is as in Lemma $\ref{natii}$.

  \end{lem}

\begin{proof}
 We
have: \begin{multline*} R_1 \odot R_2 = D^{(1)} D^{(2)} + D^{(1)}
C^{(2)} \odot (I-\zeta A^{(2)})^{-\odot}
 \odot \zeta B^{(2)} +\\ + C^{(1)}\odot (I-\zeta A^{(1)})^{-\odot} \odot
\zeta B^{(1)} D^{(2)}+\\ + C^{(1)} \odot (I-\zeta
A^{(1)})^{-\odot}\ \odot\ \zeta B^{(1)}C^{(2)}
 \odot (I-\zeta A^{(2)})^{-\odot}\ \odot\ \zeta B^{(2)}\\ = D^{(1)}D^{(2)}
+ \begin{pmatrix} C^{(1)} & D^{(1)}C^{(2)}
\end{pmatrix} \odot \begin{pmatrix} \alpha^{-\odot} &
-\alpha^{-\odot}\odot\beta\odot\gamma^{-\odot}\\ 0 &
\gamma^{-\odot}\end{pmatrix}\odot
 \begin{pmatrix} \zeta B^{(1)}D^{(2)}\\ \zeta B^{(2)}
\end{pmatrix}, \end{multline*} where \begin{align*} \alpha
&= I-\zeta A^{(1)},\\ \beta &= -\zeta B^{(1)}C^{(2)},\\ \gamma &=
I-\zeta A^{(2)}.  \end{align*} Using the formula \[
\begin{pmatrix}
\alpha^{-\odot} & -\alpha^{-\odot}\odot\beta\odot\gamma^{-\odot}\\
0 & \gamma^{-\odot}\end{pmatrix} = { \begin{pmatrix} \alpha &
\beta\\ 0 & \gamma \end{pmatrix}}^{-\odot}, \] we have:
\begin{multline*} R_1 \odot R_2 = D^{(1)}D^{(2)}+\\ +
\begin{pmatrix} C^{(1)} & D^{(1)}C^{(2)}\end{pmatrix}
\odot{\begin{pmatrix} I-\zeta A^{(1)} & -\zeta B^{(1)} C^{(2)}\\ 0
& I-\zeta A^{(2)} \end{pmatrix}}^{-\odot}  \odot \begin{pmatrix}
 \zeta B^{(1)}D^{(2)}\\ \zeta B^{(2)} \end{pmatrix}\\
=D^{(1)}D^{(2)}+\\ +\begin{pmatrix}
C^{(1)}&D^{(1)}C^{(2)}\end{pmatrix}\odot \left(I-\zeta U
\begin{pmatrix}A^{(1)}&B^{(1)}C^{(2)}\\0&A^{(2)} \end{pmatrix}
\right)^{-\odot} \odot \zeta
U\begin{pmatrix}B^{(1)}D^{(2)}\\B^{(2)}\end{pmatrix}.
\end{multline*} In order to obtain the second formula, it is
enough to apply the first one for
\[\begin{pmatrix}R_1&I\end{pmatrix}\odot\begin{pmatrix}I\\R_2\end{pmatrix}.\]
\end{proof}

\subsection{Equivalence between the various definitions}
\label{jean-yanne1} \begin{prop}
 Definitions $\ref{defrat3}$ and $\ref{mum}$ are equivalent. \end{prop}

\begin{proof}
 Indeed, if \(R\in\mathcal{O}_{\mathbb H}^{m\times n}\) admits the
representation \eqref{voltaire1}, where \(A_i\in\mathbb{H}^{p\times
p},\) let us denote by \(\mathcal W\) the span of columns of the
 matrix-function
\(W=C\odot(I-\zeta A)^{-\odot}\). According to Theorem
\ref{findiminv}, the  finite-dimensional space \(\mathcal W\) is
resolvent-invariant, and
 \(\forall v\in\mathbb{H}^{n}\) the functions
\[
G_k=C\odot(I-\zeta A)^{-\odot}B_kv\in\mathcal{W},\ k=1,2,3,\]
are a solution of the Gleason problem for \(Rv\).\\

Conversely, assume that
\(\mathcal{W}\subset\mathcal{O}^{m}_{\mathbb{H}}\) is a
finite-dimensional resolvent-invariant space, in which \(\forall
v\in\mathbb{H}^{n}\) the Gleason problem for \(Rv\)  is solvable.
According to Theorem \ref{findiminv}, there exists a
matrix-function of the form \(W=C\odot(I-\zeta A)^{-\odot}\),
whose columns span \(\mathcal{W}\). Hence there exist constant
matrices \(B_k\) such that
\[R-R(0)=\sum_{k=1}^3\zeta_k\odot WB_k.\]
Since the hyperholomorphic variables (and more, generally, all the Fueter
polynomials) belong to the center of the ring \(\mathcal{O}_{\mathbb H}\),
we obtain for \(R\) the representation \eqref{voltaire}
with \(D=R(0)\).
\end{proof}

\begin{prop} Definitions $\ref{th1111}$ and $\ref{mum}$ are
equivalent.  \end{prop}

\begin{proof}
First of all we note that, in view of Lemmas \ref{nato} and
\ref{nati}, the space of elements of \(\mathcal{O}_{\mathbb{H}}\)
which admit the representation \eqref{voltaire} is a  subring,
 which is closed under the C-K-inversion. Substituting in
\eqref{voltaire} \(B_k=0\) (respectively, \(C=B_k=1,\ D=A_k=0\)), we
see that this subring contains constant
functions (respectively, hyperholomorphic variables),
hence it also contains
\(\mathcal{Q}_{\mathbb{H}}\).
In other words, every function, rational in the sense of
Definition \ref{th1111}, admits the representation \eqref{voltaire}.\\

In order to prove the converse implication, it suffices to
show   that every entry of
\[\left(I-\sum_{k=1}^3\zeta_kA_k\right)^{-\odot}\] belongs to
\(\mathcal{Q}_{\mathbb{H}}\). We proceed by induction on \(\dim
A_k.\) Denote \[A_k=\begin{pmatrix}\check A_k & \check a_k\\ \hat
a_k&\hat A_k\end{pmatrix}.\] Then
\[I-\sum_{k=1}^3\zeta_kA_k=\begin{pmatrix}\alpha&\beta\\
\gamma&\delta\end{pmatrix}, \] where \begin{align*}
\alpha=I-\sum_{k=1}^3\zeta_k\check A_k,&\
\beta=-\sum_{k=1}^3\zeta_k\check a_k,\\
\gamma=-\sum_{k=1}^3\zeta_k\hat a_k,&\
\delta=I-\sum_{k=1}^3\zeta_k\hat A_k.  \end{align*}
Furthermore,
\begin{multline*}(I-\sum_{k=1}^3\zeta_kA_k)^{-\odot}=\\
=\begin{pmatrix}\left(\tilde\alpha\right)^{-\odot}&
-\left(\tilde\alpha\right)^{-\odot}\odot\beta\odot\delta^{-\odot}\\
-\delta^{-\odot}\odot\gamma\odot\left(\tilde\alpha\right)^{-\odot}&
\delta^{-\odot}+\delta^{-\odot}\odot\gamma\odot\left(\tilde\alpha
\right)^{-\odot}
\odot\beta\odot\delta^{-\odot}\end{pmatrix}, \end{multline*} where
\[\tilde\alpha=\alpha-\beta\odot\delta^{-\odot}\odot\gamma.\] By
definition and the induction assumption, all the entries of
\(\left(\tilde\alpha\right)^{-\odot},\beta,\gamma,\delta^{-\odot}\)
belong to \(\mathcal{Q}_{\mathbb{H}}\).  \end{proof}

\subsection{Rational functions of two complex variables}\label{generali}
Writing \[x=z_1+z_2{\bf e_2},\text{  where }z_1=x_0+x_1{\bf e_1},
\ z_2=x_2+x_3{\bf e_1},\] one can identify the  skew field of
quaternions \(\mathbb{H}\) with the two-dimensional complex space
\(\mathbb{C}^2,\) endowed with the special structure where, in
particular,
\[z{\bf e_2}={\bf e_2}\overline z.\] The complex variables \(z_1\)
and \(z_2\)
 have the following properties: \(z_1\)
is (both right and left) hyperholomorphic, \(z_2\) is
right-hyperholomorphic, \(\overline{z_2}\) is left-hyperholomorphic.
 It holds that
\[z_1(x)=\zeta_1(x)\mathbf{e_1},\
\overline{z_2}(x)=\zeta_2(x)-\zeta_3(x)\mathbf{e_1},\qquad x\in\mathbb R^4.\]
Moreover, it follows from  \eqref{Taco} that
\[z_1^m\overline{z_2}^n=z_1^{\odot m}\odot \overline{z_2}^{\odot n}.\]
From here we conclude that complex-valued functions of two complex
variables \(z_1\) and \(\overline{z_2}\), holomorphic in a
neighborhood of the origin, are also left-hyperholomorphic, for
which the C-K-product and the point-wise product coincide. It
follows that rational functions of \(z_1\) and \(\overline{z_2}\),
holomorphic in a neighborhood of the origin,   are also rational
in the sense of our Definitions \ref{th1111} -- \ref{defrat3}.

\section{Quaternionic Arveson space}
\label{hardy}
\subsection{Positive rational kernel}
\label{leiden}

In this section
we define and study what we believe to be the appropriate counterpart
of the Arveson space of the unit ball (see Section \ref{intro})
 in the setting of
left-hyperholomorphic  functions.

 To begin with, let us  recall  Definition \ref{bsdef} of
the backward-shift operators \(\mathcal{R}_n\) and formulate the
following

\begin{prop} The common eigenvectors of the backward-shift operators
\(\mathcal R_1\), \(\mathcal R_2\), \(\mathcal R_3\)
are  functions of the form
\[
(1-\zeta a)^{-\odot},\]
 where
\(a\in\mathbb H^3\).   \end{prop}

\begin{proof}
This is a special case of Theorem \ref{findiminv}: we are looking
for \(1\)-dimensional backward-shift-invariant spaces.
\end{proof}
Set
\begin{align}\Omega&=\left\{x\in{\mathbb H}:
3x_0^2+x_1^2+x_2^2+x_3^2<1\right\},\\ \label{michele}
k_y(x)&=(1-\zeta\zeta(y)^*)^{-\odot}(x),\qquad y\in{\Omega}.
\end{align} According to Definition \ref{th1111}, the
left-hyperholomorphic function $k_y$ is rational. In view of Lemma
\ref{kernel},  we have
\begin{equation}\label{miche}
k_y(x)=\sum_{\nu\in{\mathbb Z}^3_+}\frac{|\nu|!}{\nu!}
\zeta^\nu(x)\overline{\zeta^\nu(y)}.\end{equation} The function
$k_y(x)$ is therefore positive in $\Omega$ and there exists an
associated right-linear reproducing kernel Hilbert space which is
an extension of the classical Hardy space; see \cite{as3}. As a
direct consequence of the power expansion \eqref{miche} we obtain:

\begin{thm} The reproducing kernel right-linear Hilbert space $\mathbf{H}(k)$
 with reproducing kernel $k_y(x)$ (we shall call this space
the (left) quaternionic Arveson space) is the set of functions of the
form \eqref{uuab}
 endowed
with the ${\mathbb H}$-valued inner product
\begin{equation}\label{arinprod}
\langle
f,g\rangle=\sum_{{\nu\in{\mathbb
Z}^3_+}}\frac{\nu!}{|\nu|!}\overline{g_\nu}f_\nu.\end{equation}
 \end{thm}

\begin{rem}
We note that $\zeta^\nu$ and $\zeta^\mu$ are
orthogonal in the Arveson space when $\nu\not =\mu$. In the quaternionic
Hardy
space this orthogonality condition holds only when moreover $|\nu|\not
=|\mu|$.
\end{rem}

Let us consider the C-K-multiplication operators
\begin{equation}
\mathcal M_{\zeta_n}f=\zeta_n\odot f,\ n=1,2,3.
\end{equation}

\begin{prop}\label{mulbou}
For \(n=1,2,3\) the operator \(\mathcal M_{\zeta_n}\)
is a contraction from the space \(\mathbf H(k)\) into itself.
\end{prop}

\begin{proof}
 For \(f\in\mathbf H(k)\) we have
\begin{multline*}
\langle\mathcal M_{\zeta_n}f,\mathcal M_{\zeta_n}f\rangle=
\sum_{{\nu\in{\mathbb
Z}^3_+}}\frac{(\nu+e_n)!}{|\nu+e_n|!}|f_\nu|^2
=\sum_{{\nu\in{\mathbb
Z}^3_+}}\frac{\nu_n+1}{|\nu|+1}\frac{\nu!}{|\nu|!}|f_\nu|^2\\
\leq \sum_{{\nu\in{\mathbb
Z}^3_+}}\frac{\nu!}{|\nu|!}|f_\nu|^2=\langle f,f\rangle.
\end{multline*}
\end{proof}

 Proposition \ref{mulbou} implies, in particular, that the
 C-K-multiplication operator
  \(\mathcal M_{\zeta_n}\) is a bounded linear operator from
 \(\mathbf H(k)\) into itself, hence, according to the
quaternionic version
 of the Riesz theorem (see \cite{bds} and
\cite{MR96c:30045} for more details on quaternionic Hilbert spaces
and quaternionic adjoint operators ), it has the Hilbert adjoint
\(\mathcal M_{\zeta_n}^*:\mathbf H(k)\mapsto\mathbf H(k)\),
defined by
\[\langle \mathcal M_{\zeta_n}f,g\rangle=\langle f,\mathcal
M_{\zeta_n}^*g\rangle\quad\forall f, g\in\mathbf H(k).\] The
latter turns out to coincide with
 the backward-shift operator \(\mathcal R_n\):

 \begin{prop} \label{toto123}
 \[\mathcal M_{\zeta_n}^*
={\mathcal R_n}_{|\mathbf H(k)}.\]
\end{prop}
\begin{proof}
We have \(\forall \nu\in\mathbb Z_+^3\), \(\forall\mu\geq e_n\):
\[
\langle\mathcal R_n\zeta^\mu,\zeta^\nu\rangle=\langle\frac{\mu_n}{|\mu|}
\zeta^{\mu-e_n},\zeta^\nu\rangle=\frac{(\nu+e_n)!}{|\nu+e_n|!}
\delta_{\nu}^{\mu-e_n}
=\langle\zeta^\mu,\zeta^{\nu+e_n}\rangle=
\langle\zeta^\mu,\mathcal M_{\zeta_n}\zeta^\nu\rangle.
\]
Analogously, if \(\mu_n=0\)
then
\[\langle\mathcal R_n\zeta^\mu,\zeta^\nu\rangle=0=
\langle\zeta^\mu,\zeta^{\nu+e_n}\rangle=
\langle\zeta^\mu,\mathcal M_{\zeta_n}\zeta^\nu\rangle.\]
\end{proof}

Let us denote by $\mathcal C:\mathbf H(k)\mapsto\mathbb H$ the
operator of evaluation at the origin: $\mathcal Cf:=f(0)$. Then,
in view of Proposition \ref{toto123} and Theorem \ref{To kill a
mockingbird}   for the C-K-multiplication operator \(\mathcal
M_\zeta:\mathbf H(k)^3\mapsto \mathbf H(k)\) the following
operator identity holds true
 \begin{equation}
\label{Anacostia} I-\mathcal M_{\zeta}\mathcal M_{\zeta}^*=
\mathcal C^*\mathcal C. \end{equation}

The identity \eqref{Anacostia} is the quaternionic
counterpart of \eqref{2-novembre-2003}.
In the next section we shall use it to obtain
 the counterpart of the Blaschke
factors in the quaternionic Arveson space.

\subsection{Blaschke factors}
\begin{defn}
Let \(a\in\Omega\). We define the Blaschke factor \(B_a\in\mathbf
H(k)^{1\times 3}\) by
\begin{equation}
\label{blaf}
B_a=(1-\zeta(a)\zeta(a)^*)^{\frac{1}{2}} (1-\zeta\zeta(a)^*)^{-\odot}
\odot\left(\zeta-\zeta(a)\right)
\left(I-\zeta(a)^*\zeta(a)\right)^{-\frac{1}{2}}.  \end{equation}
\end{defn}

\begin{thm}\label{rovbla}
The C-K-multiplication operator
\[\mathcal B_a=\mathcal M_{B_a}:\mathbf H(k)^3\mapsto \mathbf H(k)\]
is a contraction, and the following operator identity holds:
 \begin{equation} \label{totoche}
I-\mathcal B_a\mathcal B_a^*=({1-\zeta(a)^*\zeta(a)}) \left(I-\mathcal M_\zeta
\mathcal M_{\zeta(a)}^*\right)^{-1} \mathcal C^*\mathcal C
\left(I-\mathcal M_\zeta
\mathcal M_{\zeta(a)}^*\right)^{-*}. \end{equation}
\end{thm}

\begin{proof}
 The
proof follows the arguments of \cite{akap1}. We first note that
the operators $$I-\mathcal M_{\zeta(a)}\mathcal
M_{\zeta(a)}^*\quad{\rm  and}\quad I-\mathcal
M_{\zeta(a)}^*\mathcal M_{\zeta(a)}$$ are self-adjoint and
strictly contractive and hence the operators $$\left(I-\mathcal
M_{\zeta(a)}\mathcal M_{\zeta(a)}^*\right)^{\pm 1/2} \text{ and }
\left(I-\mathcal M_{\zeta(a)}^*\mathcal M_{\zeta(a)}\right)^{\pm
1/2}$$ are well defined. We set $$ \mathcal H:=\begin{pmatrix}
\left(I-\mathcal M_{\zeta(a)}\mathcal M_{\zeta(a)}^*\right)^{-1/2}
&-\mathcal M_{\zeta(a)} \left(I-\mathcal M_{\zeta(a)}^*\mathcal
M_{\zeta(a)}\right)^{-1/2}\\ -\mathcal M_{\zeta(a)}^*
\left(I-\mathcal M_{\zeta(a)}\mathcal M_{\zeta(a)}^*\right)^{-1/2}
&\left(I-\mathcal M_{\zeta(a)}^*\mathcal
M_{\zeta(a)}\right)^{-1/2}
  \end{pmatrix}. $$
Then  it holds that
\[\mathcal HJ\mathcal H^*=\mathcal H^*J\mathcal H=J,\]
 where $$J=\begin{pmatrix}I_{{\mathbf
H}(k)}&0\\0&-I_{{\mathbf H}(k)^3} \end{pmatrix}.$$
($\mathcal H$
is "the Halmos extension" of $-\mathcal M_{\zeta(a)}$;
see \cite{MR90g:47003}.) Thus
\begin{equation}\label{singapour}
 \mathcal C^*\mathcal C=I-\mathcal M_\zeta \mathcal M_\zeta^*
=\begin{pmatrix}I&\mathcal M_\zeta\end{pmatrix}J \begin{pmatrix}I\\
\mathcal M_\zeta^*\end{pmatrix} =\begin{pmatrix}I&\mathcal
M_\zeta\end{pmatrix}\mathcal HJ\mathcal H^*
\begin{pmatrix}I\\ \mathcal M_\zeta^*\end{pmatrix}=\mathcal XJ\mathcal X^*,
\end{equation} where
\begin{align*}
 \mathcal X&= \begin{pmatrix}\mathcal X_1 &\mathcal X_2\end{pmatrix},\\
\mathcal X_1&= (I-\mathcal M_\zeta \mathcal M_{\zeta(a)}^*)
\left(I-\mathcal M_{\zeta(a)}\mathcal M_{\zeta(a)}^*\right)^{-1/2},\\
\mathcal X_2 &= (\mathcal M_\zeta-
\mathcal M_{\zeta(a)})\left(I-\mathcal M_{\zeta(a)}^*\mathcal
M_{\zeta(a)}\right)^{-1/2}.\end{align*}
To conclude we remark that
the operator $I-\mathcal M_{\zeta(a)}\mathcal M_{\zeta(a)}^*$
is the operator of multiplication by the positive number
$1-\zeta(a)\zeta(a)^*$ and therefore commutes with all the other
operators under consideration.  Multiplying the first and the last expressions
in the equality
\eqref{singapour} by $$
\left(I-\mathcal M_{\zeta(a)}\mathcal M_{\zeta(a)}^*\right)^{1/2}(I-\mathcal M_\zeta
\mathcal M_{\zeta(a)}^*)^{-1} $$ on the left and by its adjoint on the right we
obtain \eqref{totoche}.
\end{proof}

Theorem \ref{rovbla}
allows to get some preliminary results on
interpolation in the Arveson space. Here we have:

\begin{thm}\label{interpol}
 Let $a\in \Omega$.
Then
 \begin{equation} \label{tali}
\left\{f\in{\mathcal H}(k)\ :\  f(a)=0\right\} \subset \ran\mathcal B_a.
  \end{equation} \end{thm}

\begin{proof}
The identity \eqref{totoche} in Theorem \ref{rovbla} implies
that
\[\ran(I-\mathcal B_a\mathcal B_a^*)=\spa(k_a).\]
Hence
\[\ker(I-\mathcal B_a\mathcal B_a^*)=
\left\{f\in{\mathcal H}(k)\ :\  \langle f,k_a\rangle=f(a)=0\right\}.\]
On the other hand,
\(\ker(I-\mathcal B_a\mathcal B_a^*)\subset\ran\mathcal B_a\).
\end{proof}

We note that inequality is strict in \eqref{tali}. Indeed, the space
$\ran\mathcal B_a$ is invariant under the operators $\mathcal M_{\zeta_n}$
while the set of left-hyperholomorphic functions vanishing at $a$ is not,
if  $a\not =0$.

\bibliographystyle{plain}
\def\cprime{$'$} \def\cprime{$'$} \def\cprime{$'$}

%\bibliography{/users/faculty/math/dany/bib/all}
%\bibliography{/home/user/bib/all}
%\bibliography{/root/Desktop/dany/bib/all}
%\bibliography{C:/dany/bib/all}
%\bibliography{C:/WIN98/Desktop/dany/bib/all}

\begin{thebibliography}{10}

\bibitem{agler-carthy}
J.~Agler and J.~McCarthy.
\newblock Complete {N}evanlinna-{P}ick kernels.
\newblock {\em J. Funct. Anal.}, 175:111--124, 2000.

\bibitem{adr1}
D.~Alpay, A.~Dijksma, and J.~Rovnyak.
\newblock A theorem of {B}eurling--{L}ax type for {H}ilbert spaces of functions
  analytic in the ball.
\newblock {\em Integral Equations Operator Theory}, 47:251--274, 2003.

\bibitem{ad4}
D.~Alpay and C.~Dubi.
\newblock {Backward shift operator and finite dimensional de Branges Rovnyak
  spaces in the ball}.
\newblock {\em {Linear Algebra and Applications}}, 371:277--285, 2003.

\bibitem{akap1}
D.~Alpay and H.T. Kaptano\u{g}lu.
\newblock Some finite-dimensional backward shift-invariant subspaces in the
  ball and a related interpolation problem.
\newblock {\em Integral Equation and Operator Theory}, 42:1--21, 2002.

\bibitem{assv}
D.~Alpay, B.~Schneider, M.~Shapiro, and D.~Volok.
\newblock Fonctions rationnelles et th\'eorie de la r\'ealisation: le cas
  hyper--analytique.
\newblock {\em {C}omptes {R}endus {M}ath\'ematiques}, 336:975--980, 2003.

\bibitem{as3}
D.~Alpay and M.~Shapiro.
\newblock Reproducing kernel quaternionic {P}ontryagin spaces.
\newblock {Integral Equations and Operator Theory. To appear, 2003}.

\bibitem{as2}
D.~Alpay and M.~Shapiro.
\newblock Probl\`eme de {G}leason et interpolation pour les fonctions
  hyper--analytiques.
\newblock {\em C. R. Math. Acad. Sci. Paris}, 335, 2002.

\bibitem{as1}
D.~Alpay and M.~Shapiro.
\newblock Gleason's problem and tangential homogeneous interpolation for
  hyperholomorphic quaternionic functions.
\newblock {\em {Complex Variables}}, 48:877--894, 2003.

\bibitem{MR1758582}
W.~Arveson.
\newblock The curvature invariant of a {H}ilbert module over ${{\mathbb
  C}}[z\sb 1,\cdots,z\sb d]$.
\newblock {\em J. Reine Angew. Math.}, 522:173--236, 2000.

\bibitem{btv}
J.~Ball, T.~Trent, and V.~Vinnikov.
\newblock Interpolation and commutant lifting for multipliers on reproducing
  kernel {H}ilbert spaces.
\newblock In {\em Proceedings of Conference in honor of the 60--th birthday of
  {M.A. Kaashoek}}, volume 122 of {\em Operator {T}heory: {A}dvances and
  {A}pplications}, pages 89--138. Birkhauser, 2001.

\bibitem{bgk1}
H.~Bart, I.~Gohberg, and M.~Kaashoek.
\newblock {\em Minimal factorization of matrix and operator functions},
  volume~1 of {\em {Operator {T}heory: {A}dvances and {A}pplications}}.
\newblock Birkh{\" a}user Verlag, Basel, 1979.

\bibitem{MR1903737}
V.~Bolotnikov and L.~Rodman.
\newblock Finite dimensional backward shift invariant subspaces of {A}rveson
  spaces.
\newblock {\em Linear Algebra Appl.}, 349:265--282, 2002.

\bibitem{bds}
F.~Brackx, R.~Delanghe, and F.~Sommen.
\newblock {\em Clifford analysis}, volume~76.
\newblock Pitman research notes, 1982.

\bibitem{MR80c:47010}
S.~W. Drury.
\newblock A generalization of von {N}eumann's inequality to the complex ball.
\newblock {\em Proc. Amer. Math. Soc.}, 68(3):300--304, 1978.

\bibitem{MR90g:47003}
H.~Dym.
\newblock {\em ${J}$--contractive matrix functions, reproducing kernel
  {H}ilbert spaces and interpolation}.
\newblock Published for the Conference Board of the Mathematical Sciences,
  Washington, DC, 1989.

\bibitem{fueter}
R.~Fueter.
\newblock Die {T}heorie der reg\"ularen {F}unktionen einer quaternionen
  {V}ariablen.
\newblock In {\em Comptes rendus du congr\`es international des
  math\'ematiciens, {O}slo 1936, {T}ome {I}}, 1937.

\bibitem{MR1855072}
K.~Ga{\l}kowski.
\newblock Minimal state-space realization for a class of linear, discrete,
  n{D}, {S}{I}{S}{O} systems.
\newblock {\em Internat. J. Control}, 74(13):1279--1294, 2001.

\bibitem{sprossig}
K.~G{\"u}rlebeck and W.~Spr{\"o}ssig.
\newblock {\em Quaternionic and {C}lifford calculus for physicists and
  engineers}, volume~1 of {\em Mathematical methods in practice}.
\newblock John {W}iley and {S}ons, 1997.

\bibitem{MR87j:32012}
G.~M. Henkin.
\newblock The method of integral representations in complex analysis.
\newblock In {\em Current problems in mathematics. Fundamental directions,
  Vol.\ 7}, pages 23--124, 258. Akad. Nauk SSSR Vsesoyuz. Inst. Nauchn. i
  Tekhn. Inform., Moscow, 1985.

\bibitem{kfa}
R.E. Kalman, P.L. Falb, and M.A.K. Arbib.
\newblock {\em Topics in mathematical system theory}.
\newblock Mc Graw--Hill, New-York, 1969.

\bibitem{MR95m:30064}
H.~Malonek.
\newblock Hypercomplex differentiability and its applications.
\newblock In {\em Clifford algebras and their applications in mathematical
  physics (Deinze, 1993)}, volume~55 of {\em Fund. Theories Phys.}, pages
  141--150. Kluwer Acad. Publ., Dordrecht, 1993.

\bibitem{rudin-ball}
W.~Rudin.
\newblock {\em Function theory in the unit ball of ${\mathbb{C}}^n$}.
\newblock Springer--{V}erlag, 1980.

\bibitem{MR96c:30045}
M.~V. Shapiro and N.~L. Vasilevski.
\newblock Quaternionic $\psi$-hyperholomorphic functions, singular integral
  operators and boundary value problems. {I}. $\psi$-hyperholomorphic function
  theory.
\newblock {\em Complex Variables Theory Appl.}, 27(1):17--46, 1995.

\bibitem{MR54:11066}
E.L. Stout.
\newblock {\em The theory of uniform algebras}.
\newblock Bogden \& Quigley, Inc., Tarrytown-on-Hudson, N. Y., 1971.

\bibitem{MR80g:30031}
A.~Sudbery.
\newblock Quaternionic analysis.
\newblock {\em Math. Proc. Cambridge Philos. Soc.}, 85(2):199--224, 1979.

\bibitem{MR2001e:93002}
Eva Zerz.
\newblock {\em Topics in multidimensional linear systems theory}, volume 256 of
  {\em Lecture Notes in Control and Information Sciences}.
\newblock Springer-Verlag London Ltd., London, 2000.

\end{thebibliography}
\end{document}